\documentclass[12pt]{article}
\usepackage{amsmath,amsfonts}

\newcommand{\wt}{\mathop{\rm wt}}

\begin{document}

\centerline{\bf \LARGE Search for Three forged Coins} \vspace{.2in}
\centerline{\large Shaohui Zhai} \vspace{.1in} \centerline{szhai@lynx.neu.edu}
\centerline{Department of Mathematics}\centerline{Northeastern University}
\vspace{1em}
{\bf Abstract} {\small We give a simple and efficient sequential weighing algorithm
to search for three forged coins, with an asymptotic rate of $0.57$, which is
better than that $(0.46)$ of the known static algorithm given by Lindstr\"{o}m.
Then we construct from the weighing algorithm an zero-error code for the 
three-user Multiple Access Adder Channel with feedback with rate $0.57$.} 

\section{Problem} Let $T$ be a set of $t$ coins which are either
genuine or forged. Genuine coins have normalized weight $0$, forged coins have
normalized weight $1$, so that for any subset $W$ of $T$, the weight of $W$,
written by $\wt(W)$,
is the number of forged coins in $W$. Since one weighing can determine the number
of forged coins in $T$, we assume that the number of the forged coins in $T$
is known and denoted by $s$. 
We want to identify all the forged coins in $T$ by weighing some subsets of $T$, 
with as few weighings as possible.

In a {\it static} weighing algorithm, all the subsets chosen to weigh are 
determined in advance; in a {\it sequential} weighing algorithm, each subset
depends on all the previous weighing results.

In this paper, we give a sequential weighing algorithm for $s=3$ forged coins,
which is efficient when the total number of coins $t$ is large. For convenience,
we assume that $t=2^m$ for some $m$. If $t$ is not a power of 2, we can add 
``dummy" genuine coins to make up the number without affecting the asymptotic
efficiency of an algorithm.

\section{The Algorithm: Three Stages}\label{alg}

\noindent {\bf Stage 0.} We reduce $T$ to three
disjoint subsets, $T_1, T_2, T_3$, with exactly one forged coin in
each set and with $|T_1|=|T_2|=|T_3| =2^n$ for some $n$. First, we
bisect $T$ and weigh one of the halves. If the weighing result is $0$ or $3$,
we bisect the half with $3$ forged coins and weigh one of the quarters. 
We continue these bisections until the weighing result is 1 or 2.
Suppose $l_{1}$ steps are needed and the resulting subsets are
$T'$ and $T''$, containing one and two forged coins, respectively.
Next, we bisect $T''$ until the weighing result is
1. Suppose $l_{2}$ steps are needed and the resulting subsets are
denoted by $T_1$ and $T_2$.  Finally, we do $l_{2}$ bisections on
$T'$ and denote the resulting subset with one forged coin by
$T_3$.

Stage 0 takes $l_{1} + 2l_{2}$ weighings and $n = m-l_{1}
-l_{2}$.

Label each coin in $T_i$ by a distinct binary sequence of length $n$. Let $c$
indicate the label of a coin; $c^{(i)}$ the $i$-th digit of $c$; and
$c_1, c_2, c_3$ the labels of the three forged coins in $T_1, T_2,
T_3$ respectively. \\

\noindent {\bf Stage 1.} We do $n$ static weighings. The $i$-th weighing is
the subset $W_i$ of all coins in $T_1, T_2,$ and $T_3$ whose $i$-th digit is $1$: 
$$ W_i=\{c\in T_1\cup T_2\cup T_3\,|\, c^{(i)}=1\}. $$ 
If $\wt(W_i)=0$ (resp. $=3$), then none (resp. all) of $c_1, c_2, c_3$ are
included in $W_i$: that is, $c^{(i)}_1=c^{(i)}_2=c^{(i)}_3=0$ (resp. $=1$). If
$\wt(W_i)=1$ or $2$, then the $i$th digit of $c_1, c_2$ and $c_3$ is not
completely determined, and we call such a position an ambiguity. If $\wt(W_i)=1$,
exactly one of $c^{(i)}_1, c^{(i)}_2$, $c^{(i)}_3$ is
$1$; if $\wt(W_i)=2$, exactly one of $c^{(i)}_1, c^{(i)}_2$, $c^{(i)}_3$ is
$0$. Suppose there are $l_3$ weighings with result $1$ or $2$.
\\

\noindent {\bf Stage 2.} We use $l_3$ weighings to resolve the
$l_3$ ambiguous digits. Suppose $\wt(W_i)=1$ or $2$. To resolve
the $i$-th digits of $c_1, c_2,c_3$, we weigh the subset
$$W'_i=\{c\in T_1\,|\,c^{(i)}=1\}\cup \{c\in T_2\,|\,c^{(i)}=0\}.$$

Suppose $\wt(W_i)=1$. Then
exactly one of $c^{(i)}_1, c^{(i)}_2$, $c^{(i)}_3$
is $1$. If $\wt(W'_i)=0$, then $c_1,c_2\notin W'_i$,
thus $c^{(i)}_1=0$, $c^{(i)}_2=1$, so
$c^{(i)}_3$ has to be $0$. If $\wt(W'_i)=1$,
then either $c_1$ or $c_2$ is in $W'_i$,
that is, either $c^{(i)}_1=c^{(i)}_2=0$ or
$c^{(i)}_1=c^{(i)}_2=1$. But the latter case is
impossible, thus $c^{(i)}_1=c^{(i)}_2=0$,
so $c^{(i)}_3$ must be $1$. If $\wt(W'_i)=2$,
then $c_1,c_2\in W'_i$, thus $c^{(i)}_1=1$,
$c^{(i)}_2=0$, so $c^{(i)}_3$ must be $0$.

Suppose $\wt(W_i)=2$.  Then exactly one of $c^{(i)}_1,
c^{(i)}_2$, $c^{(i)}_3$ is $0$.  As above,
$c^{(i)}_1, c^{(i)}_2$, $c^{(i)}_3$ are
$(0,1,1), (1,1,0)$ or $(1,0,1)$,
when $\wt(W'_i)=0,1$ or $2$ respectively.

The results of Stage 1 and Stage 2 are illustrated in
the following table.\\[1em]
$$
\begin{array}{|c|c||ccc|} \hline
\;\wt(W_i)\;&\;\wt(W'_{i})\;&\;c^{(i)}_1\;&\;c^{(i)}_2 \;&\;c^{(i)}_3\;\\ \hline
 0       &             & 0         & 0         & 0 \\ \hline
 3       &             & 1         & 1         & 1 \\ \hline
 1       & 0           & 0         & 1         & 0 \\
         & 1           & 0         & 0         & 1 \\
         & 2           & 1         & 0         & 0 \\ \hline
 2       & 0           & 0         & 1         & 1 \\
         & 1           & 1         & 1         & 0 \\
         & 2           & 1         & 0         & 1 \\
\hline \end{array}
$$\\[.5em]

\section{Mean Duration}

{\bf Lemma.} The mean number $N$ of weighings in the algorithm is
$$
\sum_{l_1=1}^{m-1}\!\!
(\tfrac {1}{4})^{l_1-1} \tfrac {3}{4}
\!\left( \!l_1 \!+
      \!\!\!\!\!\sum_{l_2=1}^{m-l_1-1}\!\!\!
           (\tfrac{1}{2})^{l_2}\!
       \!\left(\!2l_2\!+\!m\!-\!l_1\!-\!l_2\!+\!\!\!\!\!\sum_{l_3=0}^{m-l_1-l_2}\!\!
             \!\!l_{3} \tbinom {m-l_1-l_2}{l_3} (\tfrac 34)^{l_3}
             (\tfrac 14)^{m-l_1-l_2-l_3} \right)\right)
$$
\vspace{0em}

\noindent {\bf Proof.} In Stage 0, $l_1$ is the number of weighings for the 
first occurrence of a weighing result of 1 or 2 to appear. Consider the first 
occurrence of a weighing result 1 or 2 to be a geometric random variable,
then the probability to have $l_1$ weighings is $\frac {3}{4}(\frac
{1}{4})^{l_{1}-1} $, since the probability of getting weighing result 1 or 2
is $\frac {3}{4}$. Similarly, the probability to have $l_2$ weighings is
$2^{-l_2}$, since the probability of getting weighing result 1 is $\frac
12$.  In Stage 1, the number of weighings is $n=m-l_{1} -l_{2}$.
In Stage 2, the probability of having $l_3$ ambiguities is
$\binom {n}{l_{3}} (\frac {3}{4}) ^{l_{3}}(\frac
{1}{4})^{n-l_{3}}$, and $l_3$ weighings is needed to
resolve the $l_3$ ambiguities.~$\bullet$\medskip

Simplifying the formula in the Lemma, we get $$ N=\frac {7}{4}m
-\frac {1}{2} - \frac {4m+22}{4^m} -\frac {24m-45}{2^{m+1}}.$$ 
We use the asymptotic rate $R$ to measure the efficiency of an algorithm when
$t$ goes to infinity. Then the rate of this weighing algorithm is 
$$ R \stackrel {\Delta}{=} \lim_{t
\rightarrow \infty}\frac{\log_2 t}{N}=\lim_{m \rightarrow \infty}
\frac{m}{N}=\frac{4}{7}=0.571429.$$ 

The best known weighing algorithm is Lindstr\"{o}m's \cite{Lind} static algorithm
which gives $R=0.46$.
The theoretical rate bound for static weighing algorithm is 
$0.6$ \cite{Dyachkov}. The
theoretical rate bound for sequential weighing algorithm is unknown.

\section{Corresponding code for three-user
Adder Channel with noiseless feedback}

The three-user Multiple Access Adder Channel takes three independent binary 
inputs and outputs a ternary symbol which is the sum of the three inputs.                
The communication system of this channel with noiseless feedback 
is illustrated in Fig. 1.\\[0.2in]

\unitlength 0.90mm \linethickness{0.4pt}
\begin{picture}(154.01,56.66)
\put(19.00,30){\framebox(20.00,5.67)[cc]{\small Encoder  $3$}}
\put(19.00,40.33){\framebox(20.00,5.67)[cc]{\small Encoder  $2$}}
\put(19.00,50.66){\framebox(20.00,5.67)[cc]{\small Encoder  $1$}}
\put(72.00,30){\framebox(12,25)[cc] {$+$}}
\put(97.35,35.33){\framebox(18.00,5.67)[cc]{\small Decoder}}
\put(-1.00,50.66){\framebox(11.33,5.67)[cc]{\small$[M_1]$}}
\put(-1.00,40.33){\framebox(11.33,5.67)[cc]{\small$[M_2]$}}
\put(-1.00,30){\framebox(11.33,5.67)[cc]{\small$[M_3]$}}
\put(55,55.00){\makebox(0,0)[cc]{\tiny${x^{(k)}_1(m_1,y_1,\ldots,y_{k-1})}$}}
\put(55,45.66){\makebox(0,0)[cc]{\tiny$x^{(k)}_2(m_2,y_1,\ldots,y_{k-1})$}}
\put(55,35.00){\makebox(0,0)[cc]{\tiny$x^{(k)}_3(m_3,y_1,\ldots,y_{k-1})$}}
\put(90.01,40.83){\makebox(0,0)[cc] {\tiny$y_k$}}
\put(132.01,38.00){\makebox(0,0)[cc] {\tiny$d(y_1,y_2,\ldots)$}}
\put(10.34,52.70){\vector(1,0){8.33}}
\put(10.34,43.20){\vector(1,0){8.33}}
\put(10.34,32.00){\vector(1,0){8.33}}
\put(39.34,52.70){\vector(1,0){32.33}}
\put(39.00,43.20){\vector(1,0){32.67}}
\put(39.34,32.33){\vector(1,0){32.33}}
\put(115.68,38.00){\vector(1,0){5.67}}
\put(79.34,5.67){\makebox(0,0)[cc] {{\small Fig. 1. \,Three-user Multiple
Access Adder Channel with feedback}}}
\put(84.67,38.00){\vector(1,0){12.33}}
\put(90.33,37.67){\vector(0,-1){12.00}}
\put(90.33,25.67){\line(0,-1){4}}
\put(60,24){\makebox(0,0)[cc]{\tiny noiseless feedback link}}
\put(13.33,45.33){\makebox(0,0)[lc]{\tiny $m_2$}}
\put(13.33,35.67){\makebox(0,0)[lc]{\tiny $m_3$}}
\put(13.00,54.33){\makebox(0,0)[lc]{\tiny $m_1$}}
\put(90.33,21.67){\line(-1,0){75.00}}
\put(15.33,21.67){\line(0,1){29.33}}
\put(15.33,51){\vector(1,0){3.67}}
\put(15.33,41.00){\vector(1,0){3.33}}
\put(15.33,31.00){\vector(1,0){3.67}}
\end{picture}\\[0.1in]

In this communication, there are three message sets $[M_1]=\{1,\ldots,M_1\},$\\ 
$[M_2]=\{1,\ldots,M_2\},$ and $[M_3]=\{1,\ldots,M_3\}$. To send messages
$m_1\in [M_1], m_2\in [M_2],$ and $m_3\in [M_3]$, the encoders send
their binary codewords symbol by symbol simultaneously.
The output of each transmission is the sum of all the three input symbols. With
the noiseless feedback channel, each encoder knows all the previous output 
symbols, and then decides its next input symbol, depending on the message itself 
and all the previous outputs.

We derive a simple and efficient zero-error
code from the algorithm in Section \ref{alg} for this channel.

Let $M_1=M_2=M_3=2^l$, and $m_1\in [M_1], m_2\in [M_2],$ and $m_3\in [M_3]$ 
be sent by the three users, respectively. 
We first give each message in $[M_i]$ a unique $l$-binary sequence, denoted
by $b$. Let $b^{(k)}$ be the $k$-th digit
of $b$, and $b_1, b_2, b_3$ be the $l$-binary sequences of $m_1,m_2,m_3$, 
respectively.

The code consists of two stages.\medskip

{\bf Stage 1} consists of $l$ transmissions. The three encoders send 
$b_1, b_2,b_3$ simultaneously, symbol by symbol.
The outputs of the $l$ transmissions are denoted by $y_1,\ldots,y_l\in 
\{0,1,2,3\}$. 

For $1\leq k\leq l$,
if $y_k=0$ (resp. 3), then, $b_i^{(k)}=0$ (resp. 1) for all $i$;
if $y_k=1$ or 2, then $b^{(k)}_1, b^{(k)}_2,b^{(k)}_3$ can not be completely 
determined, since we only know
that exactly one of them is 1 when
$y_k=1$, and exactly one of them is 0 when
$y_k=2$, such a position is called an ambiguity. Suppose that there 
are $a$ such ambiguities after the first
$l$ transmissions. In Stage 2, we will resolve these $a$ ambiguous positions.
\medskip

{\bf Stage 2} consists of $a$ transmissions. Let $y_k=1,$ or 2. To resolve
the $k$-th ambiguous position, the three encoders send $b^{(k)}_1, 
1-b^{(k)}_2,$ and 0, respectively. 

Let $y'_k$ denote the output of this transmission. 

If $y'_k=0,$ then the three inputs are $0,0,$ and 0. So
$b^{(k)}_1=0, b^{(k)}_2=1,$ and $b^{(k)}_3=0$ when $y_k=1$, 
$b^{(k)}_1=0, b^{(k)}_2=1,$ and $b^{(k)}_3=1$ when $y_k=2$.

If $y'_k=1,$ then either the three inputs are $1,0,0$ or $0,1,0$. When
$y_k=1$, the three inputs can only be $0,1,0$, so  
$b^{(k)}_1=0, b^{(k)}_2=0,$ and $b^{(k)}_3=1$; when
$y_k=2$, the three inputs can only be $1,0,0$, so  
$b^{(k)}_1=1, b^{(k)}_2=1,$ and $b^{(k)}_3=0$. 

If $y'_k=2,$ then the three inputs are $1,1,$ and 0. So
$b^{(k)}_1=1, b^{(k)}_2=0,$ and $b^{(k)}_3=0$ when $y_k=1$; 
$b^{(k)}_1=1, b^{(k)}_2=0,$ and $b^{(k)}_3=1$ when $y_k=2$.\medskip

The decoding is illustrated in the following table.\\[1em]
$$
\begin{array}{|c|c||ccc|} \hline
\;\;y_k\;\;&\;\;y'_k\;\;&\;\;b^{(k)}_1\;\;&\;\;b^{(k)}_2 \;\;&\;\;b^{(k)}_3\;\;\\ \hline
 0       &             & 0         & 0         & 0 \\ \hline
 3       &             & 1         & 1         & 1 \\ \hline
 1       & 0           & 0         & 1         & 0 \\
         & 1           & 0         & 0         & 1 \\
         & 2           & 1         & 0         & 0 \\ \hline
 2       & 0           & 0         & 1         & 1 \\
         & 1           & 1         & 1         & 0 \\
         & 2           & 1         & 0         & 1 \\
\hline \end{array}
$$\\[.5em]

The asymptotic transmission rate for each user is the same:
$$R_i=\lim_{M_i\to \infty}\frac{\log_2M_i}{\mbox{average number of transmissions}}
=\lim_{l\to\infty}\frac{l}{l+\frac{3l}{4}} = \frac 47 = 0.571429.$$

\end{document}